\documentclass{article}

\usepackage{amsmath}
\usepackage{amssymb}
\usepackage{epsfig}
\usepackage{tabularx}
\usepackage{url}

\newcommand{\x}{\textbf{x}}
\renewcommand{\v}{\textbf{v}}
\renewcommand{\k}{\textbf{k}}
\renewcommand{\d}{\partial}
\newcommand{\R}{\mathbb{R}}
\renewcommand{\u}{\textbf{u}}
\newcommand{\eps}{\varepsilon}
\newcommand{\sinc}{\text{sinc}}
\newcommand{\Qt}{\tilde{Q}}
\newcommand{\dt}{\Delta t}
\newcommand{\dx}{\Delta x}

\begin{document}

\title{A hybrid OpenMP and MPI implementation of a conservative spectral method for the Boltzmann equation}

\author{Jeffrey R. Haack \thanks{Department of Mathematics, The University of Texas at Austin, 2515 Speedway, Stop C1200 Austin, Texas 78712}}

\maketitle

\begin{abstract}
We demonstrate the implementation of a hybrid OpenMP and MPI parallelization of a conservative spectral method for the Boltzmann equation originally developed by Gamba and Tharkabhushaman. We perform a scaling analysis to demonstrate that the problem is well suited to parallelization, and find that the computational time scales linearly with the number of compute nodes on high performance computing resources. The original method has also been improved to higher order in space and time and is implemented on non-uniform grids in physical space. We test this scheme for an example problem in which a kinetic boundary layer generates a shock wave for large space and long times. This is the first time that the fully nonlinear Boltzmann collision operator has been used to compute this problem.
\end{abstract}


\section{Introduction} \label{sec:intro}

In this paper we consider the solution of the fully nonlinear space inhomogenous Boltzmann equation, which is an integrodifferential equation that describes the evolution of a dilute gas where only binary interactions between particles are considered. In particular, we extend the conservative spectral method developed by Gamba and Tharkabhushanam \cite{GamTha09, GamTha10} for use on high performance computing resources using OpenMP and MPI, and study the evolution of a kinetic boundary layer problem that was originally computed by Aoki et al \cite{AokiSoneNisSug91} using the Bhatnagar-Gross-Krook (BGK) approximation for collisions.  There are many difficulties associated with numerically solving the Boltzmann equation, most notably the dimensionality of the problem and the conservation of the collision invariants. For physically relevant
three dimensional applications the distribution function is seven dimensional and the velocity domain is unbounded. In addition, the
collision operator is nonlinear and requires evaluation of a five dimensional integral at each point in phase space. The collision
operator also locally conserves mass, momentum, and energy, and any
approximation must maintains this property to ensure that macroscopic
quantities evolve correctly.

Most numerical computation of the collision operator is based on
stochastic Monte Carlo methods, for example the methods of Bird
\cite{Bird} and Nanbu \cite{Nanbu}. These methods strongly reduce the
dimensionality of the problem by approximating the collision mechanism
through a stochastic, particle-based description, which ensures that
time and memory is not wasted on computing regions where $f$ is near
zero. Because these methods use the microscopic collision mechanism
directly they exactly conserve the quantities such as mass or
energy. However, the stochastic nature of Monte Carlo methods give
rise to statistical fluctuations in the numerical results. The amount of computational resources required grow very quickly in nonsteady problems, as well as problems where  collisions strongly dominate the system, flows with high mean velocities, and
problems with a distribution function that is far from equilibrium.
In recent years deterministic methods have been able to catch up with Monte Carlo
solvers, obtaining similar results with relatively coarse meshes and
providing a faster rate of convergence without having to worry about
the fluctutations or how close they are to Maxwellians
\cite{CheGamMajShu10}.

To avoid these problems, a deterministic class of so called discrete
velocity methods (DVM) were proposed. Rather than using a random
collection of particles to collide, these methods discretize the
velocity space on a finite grid in such a way that the collision
mechanics are satisfied exactly for pairs of velocity grid
points. This method originated with Broadwell \cite{Broadwell64}, and
was later extended by others \cite{Buet96, BroGolStu89, InaStu90,
  RogSch94}. This approach requires careful choice of grid points in
velocity space and on the unit sphere to preserve the conserved
quantities. However, while these methods are able to satisfy the
conservation properties of the collision operator, many have observed
accuracy of only $O(N^{-d/2})$, where $N$ is the total number of
velocity grid points in each dimension, despite being as
computationally expensive as naively applying a standard quadrature
formula directly to the collisional integral without any regard for
conservation \cite{BobPalSch97, PanHei02}. Recent work by the group of Varghese \cite{VargheseDVM} uses a combined Monte-Carlo and discrete velocity formulation, which selects random pairs of velocity grid points for collision, then conservatively interpolates the result back onto the grid. This method is less noisy than classical DSMC, and outperforms it in regions where there is less activity.

Spectral methods are a deterministic approach that compute the collision operator to high accuracy by exploiting its Fourier
structure. These methods grew from the analytical works of Bobylev \cite{Bob88} developed for the Boltzmann equation for 
 Maxwell type potential interactions and integrable angular cross section, where the corresponding Fourier transformed equation takes a 
closed form. Spectral methods provide many advantages over Direct Simulation Monte Carlo Methods (DSMC) because they are more suited to time dependent problems, low Mach number flows, high mean velocity flows, and flows that are away from equilibrium. In addition, deterministic methods avoid the statistical fluctuations that are typical of particle based methods.  Spectral approximations for this type of models where first proposed by Pareschi and Perthame \cite{ParPer96}.
Later Pareschi and Russo \cite{ParRus00} applied this work to variable hard potentials by periodizing the problem and its solution and implementing spectral collocation methods.
 
Using this representation, Bobylev and Rjasanow developed
methods for the case of Maxwell molecules \cite{BobRja97} and hard
spheres \cite{BobRja99} to derive the convolution and approximate it
with numerical quadrature. Ibragimov and Rjasanow extended these ideas
to the more general case of variable hard spheres \cite{IbrRja02} and
simplified the integration domain by explicitly computing the
spherical integral in the derivation of the weight function for the
convolution.  In a related work, Pareschi and Russo \cite{ParRus00}
applied this framework to develop methods for the Maxwell molecules,
hard spheres, and variable hard spheres cases using a collocation
method, which uses orthogonal polynomials to reduce the convolution
integral to a convolution sum, providing spectral accuracy to the
collision integral computation.

Inspired by the work of Ibragimov and Rjasanow \cite{IbrRja02}, Gamba and Tharkabhushanam \cite{GamTha09, GamTha10} observed that the Fourier transformed collision operator takes a simple form of a weighted convolution and  developed a spectral method based on the weak form of the Boltzmann equation that provides a general framework for computing both elastic and inelastic collisions. Macroscopic conservation is enforced by solving a numerical constrained optimization problem that finds the closest distribution function in $L_2$ to the output of the collision term that conserves the macroscopic quantities.  These methods do not rely on periodization but rather on the use of the fast Fourier transform in the computational domain, and convergence to the solution of the continuous problem is obtained by the use of the extension operator in Sobolev spaces \cite{AlonGamTha}.   

All of these methods are computationally expensive, so there is a natural desire to extend them to high performance computing resources. However, there have been relatively few published works on parallelization of Boltzmann solvers. Graphics Processing Units (GPUs) have emerged as a resource for applications other than graphics, providing hundreds to thousands of lightweight independent processors already geared to independently operating on large sets of data. Frezzotti and collaborators \cite{FrezGhiGib11,FrezGhiGib11-2,FrezGhiGib11-3} implemented DSMC and BGK type Boltzmann models on GPUs, and in parallel Malkov and Ivanov and collaborators have explored implementing classical DSMC as well as the discrete velocity Monte Carlo of Varghese et al. on GPUs  \cite{MalIvanovGPU,MalPolIva12,KasSheIva12}. Many existing codes cannot be ported directly to GPUs, due to their relatively simple instruction set and complicated memory management structure. In 2012, Intel released a new architecture known as Many Integrated Core (MIC). This architecture combines many more powerful cores on a single chip than ever before, all able to access the same pooled memory. These cores are classical CPUs, rather than GPUs, which allows relatively straightforward speedup of existing codes \cite{GamHaaRGD}.

This paper presents the extension of the deterministic spectral method of Gamba and Tharkabhushanam to CPU-based parallelization using the OpenMP and MPI APIs, as opposed to GPUs. We also present a second order in time and space extension of this method that implements nonuniform grids. Sections 2 and 3 are presented for the benefit of the reader; much of the spectral formulation can be found in \cite{GamTha09,GamTha10}. This method has been tested on the Ranger supercomputer at the Texas Advanced Computing Center (TACC) and timing studies are provided to explore the scaling of the method to more and more processors for large problems, and will serve as a stepping stone to future implementation on MIC architecture. We present some $1D$ in physical space results for a problem proposed by Aoki et al.\cite{AokiSoneNisSug91} where a sudden change in wall temperature results in a shock that cannot be explained by classical hydrodynamics. 
\section{The space inhomogeneous Boltzmann equation} \label{sec:BTE}
The space inhomogeneous initial-boundary value problem for the Boltzmann equation is given by
\begin{equation}\label{BTE}
\frac{d}{dt} f(\x,\v,t) + \v \cdot \nabla_x f(\x, \v, t) = \frac{1}{\eps} Q(f,f),
\end{equation}
with
\begin{align*}
&\x \in \Omega \subset \R^d, \qquad \v \in \R^d \\
&f(\x,\v,0) = f_0(\x,\v) \\
&f(\x,\v,t) = f_B(\x,\v,t), \qquad \x \ \in \d \Omega.
\end{align*}
where $f(\v,t)$ is a probability density distribution in $\v$-space and $f_0$ is assumed to be locally integrable with respect to $\v$ and the spatial boundary condition $f_B$ will be specified in below. The dimensionless  parameter $\eps > 0$ is the scaled Knudsen number, which is defined as the ratio between the mean free path between collisions and a reference macroscopic length scale.
 
The collision operator $Q(f,f)(\x,\v,t)$ is a bilinear integral form local in $t$ and $\x$ and is given by
\begin{equation}\label{Q_general}
Q(f,f)(\cdot,\v,\cdot) = \int_{\v_\ast \in \R^d} \int_{\sigma \in S^{d-1}} B(|\v - \v_\ast|,\cos \theta) (f(\v_\ast ')f(\v') - f(\v_\ast)f(\v)) d\sigma d\v_\ast,
\end{equation}
where the velocities $\v', \v_\ast'$ are determined through a collision rule \eqref{velocity_interact}, depending on $\v, \v_\ast,$, and the positive term of  the integral in \eqref{Q_general} evaluates  $f$ in the pre-collisional velocities that will take the direction $\v$ after an interaction.  The collision kernel $B(|\v - \v_\ast|,\cos \theta)$ is a given non-negative function depending on the size of the relative velocity $\u := \v - \v_\ast$ and $\cos \theta = \frac{\u \cdot \sigma}{|\u|}$, where $\sigma$ in the $n-1$ dimensional sphere $S^{n-1}$ is referred as the scattering direction  
of the post-collisional elastic relative velocity.

For the following we will use the velocity elastic (or reversible) interaction law in center of mass-relative velocity coordinates
\begin{align} 
&\v' = \v + \frac{1}{2}(|\u|\sigma - \u), \qquad \v_\ast ' = \v_\ast - \frac{1}{2}(|\u|\sigma - \u) \label{velocity_interact} \\
&B(|\u|,\cos \theta) = |\u|^\lambda b(\cos \theta) \, \nonumber .
\end{align}
We assume that the differential cross section function $b(\cos \theta)$ is integrable with respect to $\sigma$ on $S^{d-1}$, referred to as the Grad cut-off assumption, and that it is renormalized such that 
\begin{equation} \label{cut-off}
\int_{S^{d-1}} b(\cos \theta) d\sigma = 1.
\end{equation}
The parameter $\lambda$ regulates the collision frequency as a function of the relative speed $|\u|$. This corresponds to the interparticle potentials used in the derivation of the collisional kernel and are referred to as variable hard potentials (VHP) for $0 < \lambda < 1$, hard spheres (HS) for $\lambda = 1$, Maxwell molecules (MM) for $\lambda = 0$, and variable soft potentials (VSP) for $-3 < \lambda < 0$. The $\lambda = -3$ case corresponds to a Coulombic interaction potential between particles. If $b(\cos \theta)$ is independent of $\sigma$ we call the interactions isotropic (like the case of hard spheres in three dimensions). 

Depending on the nature of the collisions, the collision operator can have a number of collision invariants. In the case of the classical Boltzmann equation with elastic collisions, according to the Boltzmann Theorem the only collisional invariants are polynomials of the form $A + B \cdot \v + C |\v|^2$. These give rise to the classical macroscopic conserved quantities
\begin{align} \label{moments}
 \rho(\x,t) &= \int_{\v} f(\x,\v,t) d\v \qquad \text{(density)}\nonumber  \\
 \rho(\x,t) \textbf{V}(\x,t) &=  \int_{\v} \v f(\x,\v,t) d\v \qquad \text{(momentum)} \nonumber \\
 \rho(\x,t) e(\x,t) &= \frac{1}{2} \int_{\v} |\v|^2 f(\x,\v, t) d \v \qquad \text{(kinetic energy density)} 
\end{align}
%

\subsection{Boundary conditions} \label{sec:BC}
On the spatial boundary $\d \Omega$ we use a diffusive Maxwell boundary condition which is given by, for $\x \in \d \Omega$,
\begin{align} \label{BC}
f(\x,\v,t)  &= \frac{\sigma_w}{(2\pi R T_w)^(d/2)} \exp \left( - \frac{|\v - \textbf{V}_w|^2}{2 R T_w}\right), \qquad (\v - \textbf{V}_w)\cdot \textbf{n} > 0 \\
\sigma_w &= - \left(\frac{2\pi}{R T_w}\right)^{1/2} \int_{(\v - \textbf{V}_w)\cdot\textbf{n} < 0} (\v - \textbf{V}_w)\cdot\textbf{n} f(\x,\v,t) d \v, \nonumber
\end{align}
where $\textbf{V}_w$ and $T_w$ are the wall velocity and temperature, respectively, and $\textbf{n}$ is the unit normal vector to the boundary, directed into $\Omega$. The term $\sigma_w$ accounts for the amount of particles leaving the domain and ensures mass conservation in $\Omega$.


\subsection{Spectral formulation} \label{sec:spectral_cont}
The key step our formulation of the spectral numerical method is the use of the weak form of the Boltzmann collision operator. For a suitably smooth test function $\phi(\v)$ the weak form of the collision integral is given by
\begin{equation} \label{collision_weakform}
\int_{\R^d} Q(f,f) \phi(\v) d\v = \int_{\R^d \times \R^d \times S^{d-1}} f(\v)f(\v_\ast) B(|\u|,\cos \theta) (\phi(\v') - \phi(\v)) d\sigma d\v_\ast d\v
\end{equation}
If one chooses 
\begin{equation*}
\phi(\v) = e^{-i \zeta \cdot \v} / (\sqrt{2\pi})^d,
\end{equation*}
then \eqref{collision_weakform} is the Fourier transform of the collision integral with Fourier variable $\zeta$:
\begin{align}
\widehat{Q}(\zeta) &= \frac{1}{(\sqrt{2\pi})^d} \int_{\R^d} Q(f,f) e^{-i \zeta \cdot \v} d\v \nonumber \\
&= \int_{\R^d \times \R^d \times S^{d-1}} f(\v)f(\v_\ast) \frac{B(|\u|,\cos \theta)}{(\sqrt{2\pi})^d} (e^{-i \zeta \cdot \v'} - e^{-i \zeta \cdot \v}) d\sigma d\v_\ast d\v \nonumber \\
&= \int_{\R^d} G(\u,\zeta) \mathcal{F}[f(\v)f(\v-\u)](\zeta) d\u, \label{First_FT}
\end{align}
where $\widehat{[\cdot]} = \mathcal{F}(\cdot)$ denotes the Fourier transform and 
\begin{equation} \label{G_eqn}
G(\u,\zeta) = |\u|^\lambda \int_{S^{d-1}} b(\cos \theta) \left(e^{-i\frac{\beta}{2} \zeta \cdot |\u|\sigma}e^{i\frac{\beta}{2} \zeta \cdot \u} - 1\right) d\sigma
\end{equation}
Further simplification can be made by writing the Fourier transform inside the integral as a convolution of Fourier transforms:
\begin{align} \label{ContSpectral}
\widehat{Q}(\zeta) 
&= \int_{\R^d} \widehat{G}(\xi,\zeta) \hat{f}(\zeta - \xi) \hat{f}(\xi) d\xi,
\end{align}
where the convolution weights $\widehat{G}(\xi,\zeta)$ are given by
\begin{align} \label{Ghat_eqn}
\widehat{G}(\xi,\zeta) &= \frac{1}{(\sqrt{2\pi})^d}  \int_{\R^d} G(\u,\zeta) e^{-i \xi \cdot u} d\u 
\end{align}
These convolution weights can be precomputed once to high accuracy and stored for future use. For many collision types the complexity of the integrals in the weight functions can be reduced dramatically through analytical techniques. In this paper we will only consider isotropic scattering in dimension 3  $(b(\cos \theta) = 1/ 4\pi)$. In this case we have that 

\begin{align} \label{G_sinc}
G(\u,\zeta) &= \frac{|\u|^\lambda}{4\pi}\int_{S^2} e^{-i\frac{\beta}{2} \zeta |\u| \cdot \sigma} e^{i\frac{\beta}{2} \zeta \cdot \u} - 1\ d\sigma \nonumber \\
&= |\u|^\lambda \left(e^{i\frac{\beta}{2} \zeta \cdot \u} \sinc\left(\frac{\beta|\u| |\zeta|}{2}\right) - 1\right) \
\end{align}
To calculate $\widehat{G}(\xi,\zeta)$, we have
\begin{align} \label{GHat_calc}
\widehat{G}(\xi,\zeta) &= \frac{1}{(\sqrt{2\pi})^3} \int_{\R^3} G(\u,\zeta) e^{-i\xi\cdot\u} d\u \nonumber \\
&= \frac{4\pi}{(\sqrt{2\pi})^3} \int_{\R^+} r^{\lambda + 2} \left( \sinc\left(\frac{\beta r |\zeta|}{2}\right) \sinc(r|\xi - \frac{\beta}{2}\zeta|) - \sinc(r|\xi|)\right) dr.
\end{align}
%

This integral will be cut off at a point $r = r_0$, which will be determined below. Given this cutoff point, we can explicitly compute $\widehat{G}$ for integer values of $\lambda$. 
\begin{eqnarray}
\widehat{G}(\xi,\zeta) &= \frac{4\pi}{(\sqrt{2\pi})^3}&\Big(\frac{q^2(pr_0\sin(pr_0) + \cos(pr_0) - 1) - p^2(qr_0\sin(qr_0) + \cos(qr_0) - 1)}{\beta|\zeta||\xi - \frac{\beta}{2}\zeta|p^2q^2} \nonumber \\
&&- \frac{(2 - r_0^2|\xi|^2)\cos(r_0|\xi|) + 2r_0|\xi|\sin(r_0|\xi|) - 2}{|\xi|^4}\Big), \qquad \qquad \lambda = 1, \nonumber \\
\widehat{G}(\xi,\zeta) &=  \frac{4\pi}{(\sqrt{2\pi})^3}&\Big(\frac{q\sin(pr_0) - p\sin(qr_0)}{\beta|\zeta||\xi - \frac{\beta}{2}\zeta|pq}- \frac{\sin(|\xi|r_0) - |\xi|r_0\cos(|\xi|r_0)}{|\xi|^3}\Big), \ \ \ \lambda = 0. \nonumber 
\end{eqnarray}

For other values of $\lambda$, this is simply a one-dimensional integral that can be precomputed to high accuracy using numerical quadrature. The entirety of the collisional model being used is encoded in the weights, which gives the algorithm a large degree of flexibility in implementing different models. 


\section{The Conservative Numerical Method} \label{sec:numerics_setup}

\subsection{Temporal and velocity space discretization} \label{sec:discretization}

We use an operator splitting method to separate the mechanisms of collisions and advection. The system is split into the subproblems
\begin{align}
\frac{\d}{\d t} f + v \cdot \nabla_\x f = 0 \\
\frac{\d}{\d t} f = Q(f,f),
\end{align}
which are solved separately. 

Each system is evolved in time using a second-order Runge-Kutta method, and the systems are combined using Strang splitting.

In order to compute the Boltzmann equation we must work on a bounded velocity space, rather than all of $\R^d$. However typical distributions are supported on the entire domain, for example the Maxwellian equilibrium distribution. Even if one begins with a compactly supported initial distribution, each evaluation of the collision operator spreads the support by a factor of $\sqrt{2}$, thus we must use a working definition of an {\em effective support} of size  $R$ for the distribution function. Bobylev and Rjasanow \cite{BobRja99} suggested using the temperature of the distribution function, which typically decreases as $\text{exp}(-|v|^2 / 2T)$ for large $|v|$, and used a rough estimate of $R \approx 2\sqrt{2}T$ to determine the cutoff. 
We assume that the distribution function is negligible outside of a ball 
\begin{equation} \label{Ball_domain_v}
B_{R_x}(\textbf{V}(\x)) = \{\v \in \R^d : |\v - \textbf{V(\x)}| \le R_x \},
\end{equation}
where $\textbf{V}(\x)$ is the local flow velocity which depends in the spatial variable $\x$. For ease of notation in the following we will work with a ball centered at $0$ and choose a length $R$ large enough that $B_{R_x}(\textbf{V}(\x)) \subset B_R(0)$ for all $\x$.

With this assumed support for the distribution $f$, the integrals in \eqref{ContSpectral} will only be nonzero for $\u \in B_{2R}(0)$. Therefore, we set $L=2R$ and define the cube
\begin{equation} \label{Cube_domain_v}
C_L = \{ \v \in \R^d : |v_j| \le L,\,\, j = 1,\dots,d\}
\end{equation}
to be the domain of computation. For such domain, the computation of the weight function integral \eqref{GHat_calc} is cut off at $r_0=L$.

Let $N \in \mathbb{N}$ be the number of points in velocity space in each dimension. Then the uniform velocity mesh size is $\Delta v = \frac{2L}{N}$ and due to the formulation of the discrete Fourier transform the corresponding Fourier space mesh size is given by $\Delta \zeta = \frac{\pi}{L}$. 

The mesh points are defined by
\begin{align} \label{meshpoints}
v_\k &= \Delta v (\k - N/2) \nonumber \\
\zeta_\k &= \Delta \zeta (\k - N/2) \\
& \k = (k_1,\dots,k_d) \in \mathbb{Z}^d,\qquad 0 \le k_j \le N-1,\,\,\, j = 1,\dots,d
\end{align}


\subsection{Collision step discretization} \label{sec:collision_disc}
Returning to the spectral formulation \eqref{ContSpectral}, the weighted convolution integral then becomes an integral over $-\frac{\pi}{\Delta v} \le \xi_j \le \frac{\pi}{\Delta v}, \, j = 1,\dots,d$. 

Similar to what was noted in \cite{IbrRja02}, we can find the cutoff for the integration variable $\u$ through the term
\[ g(\u,\v) = f(\v)f(\v - \u) \] 
that appears in the Fourier transform term in \eqref{First_FT}. As $\textrm{supp} f = B_L(0)$, we have that
\[ \textrm{supp}\, g(\u,\cdot) = B_L(0) \cap B_L(\u) \]

To simplify notation we will use one index to denote multidimensional sums with respect to an index vector $\textbf{m}$
\begin{equation*}
\sum_{\textbf{m}=0}^{N-1} = \sum_{m_1,\dots,m_d = 0}^{N-1}.
\end{equation*}

To compute $\widehat{Q}(\zeta_\k)$, we first compute the Fourier transform integral giving $\hat{f}(\zeta_k)$ via the FFT. The weighted convolution integral is approximated using the trapezoidal rule
\begin{equation}
\widehat{Q}(\zeta_\k)= \sum_{\textbf{m} = 0}^{N-1} \widehat{G}(\xi_{\textbf{m}},\zeta_\k) \hat{f}(\xi_{\textbf{m}}) \hat{f}(\zeta_\k - \xi_{\textbf{m}}) \omega_{\textbf{m}},
\end{equation}
where $\omega_\textbf{m}$ is the quadrature weight and we set $\hat{f}(\zeta_\k - \xi_{\textbf{m}}) = 0$ if $(\zeta_\k - \xi_{\textbf{m}})$ is outside of the domain of integration. We then use the inverse FFT on $\widehat{Q}$ to calculate the integral returning the result to velocity space. 

Note that in this formulation the distribution function is not periodized, as is done in the collocation approach of Pareschi and Russo \cite{ParRus00}. This is reflected in the omission of Fourier terms outside of the Fourier domain. All integrals are computed directly only using the FFT as a tool for faster computation.The convolution integral is accurate to at least the order of the quadrature. The calculations below use the trapezoid rule, but in principle Simpson's rule or some other uniform grid quadrature can be used. However, it is known that the trapezoid rule is spectrally accurate for periodic functions on periodic domains (which is the basis of spectral accuracy for the FFT), and the same arguments can apply to functions with sufficient decay at the integration boundaries \cite{Atkinson}. These accuracy considerations will be investigated in future work. The overall cost of this step is $O(N^{2d})$. 

\subsection{Discrete conservation enforcement} \label{sec:conservation}
This implementation of the collision mechanism does not conserve all of the quantities of the collision operator. To correct this fact, we formulate these conservation properties as a constrained optimization problem as proposed in \cite{GamTha09, GamTha10}. Depending on the type of collisions we can change this constraint set (for example, inelastic collisions do not preserve energy). We focus here just on the case of elastic collisions, which preserve mass, momentum, and energy. 

Let $M = N^d$ be the total number of grid points, let $\tilde{\textbf{Q}} = (\Qt_1, \dots, \Qt_M) ^T$ be the result of the spectral formulation from the previous section, written in vector form, and let $\omega_j$ be the quadrature weights over the domain in this ordering. Define the integration matrix
\begin{equation*}
\textbf{C}_{5\times M} = \left(\begin{array}{c} \omega_j \\ v_j^i \omega_j \\ |\v_j|^2 \omega_j \end{array} \right),
\end{equation*}
where $v^i,\, i=1,2,3$ refers to the $i$th component of the velocity vector. Using this notation, the conservation method can be written as a constrained optimization problem. 

\begin{equation} 
\text{Find } \textbf{Q} = (Q_1,\dots,Q_M)^T \text{ that minimizes } \frac12 \|\tilde{\textbf{Q}} - \textbf{Q}\|_2^2 \text{ such that } \textbf{C} \textbf{Q} = \textbf{0}
\end{equation}

The solution is given by
\begin{align}
\textbf{Q} &= \tilde{\textbf{Q}} + \textbf{C}(\textbf{C}\textbf{C}^T)^{-1} \textbf{C} \tilde{\textbf{Q}} \nonumber \\
&:= \textbf{P}_N \tilde{\textbf{Q}}
\end{align}

Overall the collision step in semi-discrete form is given by
\begin{equation}
\frac{\d \textbf{f}}{\d t} = \textbf{P}_N \tilde{\textbf{Q}}
\end{equation}

The overall cost of the conservation portion of the algorithm is a $O(N^d)$ matrix-vector multiply, significantly less than the computation of the weighted convolution.


\subsection{Spatial and Transport discretization} \label{sec:transport}

For simplicity this will be presented in 1D in space, though the ideas apply to higher dimensions. In this case the transport equation reduces to
\begin{equation*}
\frac{\d f}{\d t} (x,\v, t) + v_1 \frac{\d}{\d x} f(x,\v,t) = 0.
\end{equation*}

We partition the domain into cells of size $\dx_j$ (not necessarily uniform) with cell centers $x_j$. Using a finite volume approach, we integrate the transport equation over a single cell to obtain
\begin{align}
\frac{f^{n+1}_j (v_i) - f^n_j (v_i)}{\dt} + \frac{F^n_{j+1/2} - F^n_{j-1/2}}{\dx_j} = 0, \label{transport}
\end{align}
where $t^n = n\dt$ and $F^n_{j\pm1/2}$ is an approximation of the edge fluxes $v_1f$ of the cell between time $t^n$ and $t^{n+1}$. We use a second order upwind scheme defined by 
\begin{equation}\label{upwind}
F^n_{j+1/2} = \begin{cases} v_1 (f^n_j + \frac{\dx}{2}\sigma^n_j), \qquad &v_1 \ge 0 \\
		                         v_1(f^n_{j+1} - \frac{\dx}{2}\sigma^n_{j+1}), \qquad &\text{otherwise,} \end{cases}
\end{equation}
where $\sigma_j$ is a cell slope term used in the reconstruction defined by the minmod limiter.
\begin{equation} \label{minmod}
\sigma_j = \text{minmod}\left(\frac{f^n_{j+1} - f^n_j}{x_{j+1} - x_j}, \frac{f^n_{j} - f^n_{j-1}}{x_j - x_{j-1}}, \frac{f^n_{j+1} - f^n_{j-1}}{x_{j+1} - x_{j-1}}\right)
\end{equation}
For reconstructions at the boundary of the physical domain, ghost cells and extrapolation are used to determine the reconstructed slope \cite{RedLeVeque}.

On wall boundaries the incoming flux is determined using the diffusive reflection formula \eqref{BC}. For problems without meaningful boundary interactions (e.g. shocks), a no-flux boundary condition is applied for the incoming characteristics.

\section{Parallelization}

\subsection{Shared Memory Parallelization}
The most computationally expensive term to evaluate in one time step of the Boltzmann solver is the weighted convolution in $\hat{Q}$. This requires computing the sum of $N^3$ terms: $\sum_{\textbf{m} = 0}^{N-1} \widehat{G}(\xi_{\textbf{m}}, \zeta_{\textbf{k}}) \hat{f}(\xi_{\textbf{m}}) \hat{f}(\zeta_{\textbf{k}} - \xi_{\textbf{m}})\omega_\textbf{m} $ for each of the $N^3$ values of $\zeta_\textbf{k}$ on the numerical grid. Each sum draws from essentially the same data and recombines it in a different way, which makes it ideal for shared memory parallelization.

Shared memory parallelization refers to utilizing an architecture in which multiple processes can simultaneously access the same address space in memory. This type of parallelization has grown more prominent over the past decade as multicore processors are becoming more and more prevalent in typical personal computers, and they are the building blocks of larger scale high performance computing systems. The primary motivation for this type of architecture is that memory access speed is much less than the speed of a processor, and having a common pool of memory mitigates the need for memory transfers between processors of the system. In this architecture, computational work is divided into a discrete number of threads that are distributed to the available processors, which compute each thread separately. Logistically, the most difficult hurdles to overcome in this computing paradigm are race conditions, where one thread may not recieve the correct data depending on whether another thread has modified it before access, load balancing, where work is not distributed equally among the threads, and blocking, where one thread must wait on another thread to complete some work before it can begin.

We utilize the OpenMP API \cite{openMP08} to parallelize the weighted convolution. This interface is available on most computing architectures and compilers, and in many cases can give a significant speedup for only adding two or three simple lines of code. Typically the same source code can be executed on anything ranging from a personal laptop to a high performance computing node, which illustrates the great portability of this model. 

For our application, the hurdles associated with parallelization mentioned above are not a problem. We avoid load balancing problems by evenly subdividing the work among the available cores, and each convolution is expected to take the same amount of operations as any other. Each thread makes a computation based the weight $\widehat{G}$ and distribution function $\hat{f}$, and is not dependent on communication with other threads, so there is no need to worry about race conditions or thread blocking. Taking this into account the speedup would be expected to scale linearly with $p$, the number of cores that can access the memory, however memory access between the shared memory units in the node still present latency issues on the total speedup, as observed in the test below.

In Tables 1--4 we present the timing results of a single evaluation of the collision operator for $N=16$ and $N=24$ points in each velocity direction. These tests are performed on the Texas Advanced Computing Center's supercomputers Ranger and Stampede. On Ranger, a single computational node contains four AMD Opeteron quad-core 64 bit processors, for a total of 16 cores with a frequency of 2.3 GHz each, 32 GB of shared memory, and a peak performance of 128 GFLOPS per node. Stampede is TACC's newest system, officially deploying in January 2013 and uses the new Intel MIC architecture. Each computational node consists of two Intel Xeon E-5 Processors and a Xeon Phi Coprocessor (the MIC component). The E-5 processors have 8 cores each at 2.7 GHz and 32 GB of memory, for a total of 16 cores, and the Phi Coprocessor contains 61 cores at 1.1 GHz and 8GB of fast-access memory for a total of $\sim$1070 GFLOPS of performance by itself. The peak performance of the entire system was benchmarked at 3959 Teraflops in November 2012, making it the 7th fastest supercomputer in the world \cite{top500}, and has more components yet to be added. As Stampede is still in its trial period
 the necessary libraries have not been created to run this code on the coprocessor, so the following results are run only on the 16 cores of the E-5 processors.

\begin{table}[!htbp]
\begin{tabularx}{\textwidth}{XXXX}
\hline
cores & time (s) & ratio & total speedup\\
\hline
1 &   0.068 & &\\
2 &   0.044 & 1.54 & 1.54\\
4 &   0.027 & 1.63 & 2.52 \\
8 &   0.019 & 1.42 & 3.58\\
16 & 0.018 & 1.05 & 3.78\\ 
\hline
\end{tabularx}
\caption{Time for single evaluation of collision operator on one node with OpenMP. Ranger with N=16 }
\centering
\label{timesRanger16}
\end{table}

\begin{table}[!htbp]
\begin{tabularx}{\textwidth}{XXXX}
\hline
cores & time (s) & ratio & total speedup\\
\hline
1 &  0.90567&&\\
2 &  0.4899& 1.85 & 1.85\\
4 &  0.3021& 1.62 & 3.0 \\ 
8 &  0.19197& 1.57& 4.72 \\
16 &  0.1769& 1.09 & 5.12 \\
\hline
\end{tabularx}
\caption{Time for single evaluation of collision operator on one node with OpenMP. Ranger with N=24 }
\centering
\label{timesRanger24}
\end{table}

\begin{table}[!htbp]
\begin{tabularx}{\textwidth}{XXXX}
\hline
cores & time (s) & ratio & total speedup\\
\hline
1 &   0.03337&\\
2 &   0.01627 & 2.05& 2.05\\
4 &   0.00976 & 1.67& 3.41\\
8 &   0.00599 & 1.63& 5.57\\
16 & 0.00408 & 1.47& 8.18\\ 
\hline
\end{tabularx}
\caption{Time for single evaluation of collision operator on one node with OpenMP. Stampede with N=16 }
\centering
\label{timesStamp16}
\end{table}

\begin{table}[!htbp]
\begin{tabularx}{\textwidth}{XXXX}
\hline
cores & time (s) & ratio & total speedup\\
\hline
1 &   0.34386& &\\
2 &   0.17446   & 1.97& 1.97\\
4 &   0.103       & 1.69& 3.34\\
8 &   0.06107   & 1.68  & 5.63\\
16 & 0.04611  & 1.32& 7.46\\ 
\hline
\end{tabularx}
\caption{Time for single evaluation of collision operator on one node with OpenMP. Stampede with N=24 }
\centering
\label{timesStamp24}
\end{table}

\subsection{Distributed Memory Parallelization}

The total number of operations required to compute the all of the collisional terms in a single timestep is $O(MN^6)$, where $M$ is the total number of physical space grid points. When dividing the computational work, in shared memory parallelization we are restricted to the number of processors that can physically access the shared memory. Distributed memory parallelization on the other hand consists of starting multiple processes that can communicate with each other, each with their own private address space.

The main drawback of distributed memory is memory access time. Unlike shared memory, distributed memory is less local and thus requires much more time to receive information from a process that is not in its address space when compared to a similar shared memory access. Therefore, it is important to design algorithms that limit communication between processes as much as possible. Furthermore, many of the same issues from shared memory programming still apply, such as load balancing and blocking. 

Recall that the collision term in the Boltzmann equation is local in space, thus each grid point in physical space only requires $\hat{f}(\x,\cdot)$ and the precomputed weights $\widehat{G}$ to evaluate the collision term. This allows for a natural decomposition of the computational domain by separating across physical grid points. The only communication between the computational nodes is in the transport term \eqref{transport}. Each process simply sends and receives the $O(N^3)$ values of $\hat{f}(\x,\zeta)$ at the edges of the decomposed domain required to calculate the spatial flux. The communication between processes is managed by the Message Passing Interface (MPI) protocol \cite{MPI} using an interleaved ghost cells technique \cite{GroppMPI}. This mitigates the possibility of message deadlock by ensuring that all even indexed processes send or receive data at the same time, while the odd indexed processors receive data from or send data to the even processors, respectively. The code below fills the edge cells on nodes to either side. Extrapolation is used to fill ghost cells at the physical domain boundaries, which is removed from the code below for clarity of presentation.

\begin{verbatim}
  ---------------------------------
  if((rank % 2) == 0) { //EVEN NODES
     MPI_Ssend(f[nX_node+1],N*N*N,MPI_DOUBLE,rank+1,0,MPI_COMM_WORLD);

     MPI_Recv(f[1],N*N*N,MPI_DOUBLE,rank-1,0,MPI_COMM_WORLD, &status);
     
     MPI_Ssend(f[nX_node],N*N*N,MPI_DOUBLE,rank+1,0,MPI_COMM_WORLD);

     MPI_Recv(f[0],N*N*N,MPI_DOUBLE,rank-1,0,MPI_COMM_WORLD, &status);

     MPI_Ssend(f[2],N*N*N,MPI_DOUBLE,rank-1,1,MPI_COMM_WORLD);

     MPI_Recv(f[nX_node+2],N*N*N,MPI_DOUBLE,rank+1,1,MPI_COMM_WORLD, &status);

     MPI_Ssend(f[3],N*N*N,MPI_DOUBLE,rank-1,1,MPI_COMM_WORLD);

     MPI_Recv(f[nX_node+3],N*N*N,MPI_DOUBLE,rank+1,1,MPI_COMM_WORLD, &status);
  }
  else { //ODD NODES 
    MPI_Recv(f[1],N*N*N,MPI_DOUBLE,rank-1,0,MPI_COMM_WORLD, &status); 
        
    MPI_Ssend(f[nX_node+1],N*N*N,MPI_DOUBLE,rank+1,0,MPI_COMM_WORLD);

    MPI_Recv(f[0],N*N*N,MPI_DOUBLE,rank-1,0,MPI_COMM_WORLD, &status); 
   
    MPI_Ssend(f[nX_node],N*N*N,MPI_DOUBLE,rank+1,0,MPI_COMM_WORLD);
   
    MPI_Recv(f[nX_node+2],N*N*N,MPI_DOUBLE,rank+1,1,MPI_COMM_WORLD, &status);
      
    MPI_Ssend(f[2],N*N*N,MPI_DOUBLE,rank-1,1,MPI_COMM_WORLD); 
   
    MPI_Recv(f[nX_node+3],N*N*N,MPI_DOUBLE,rank+1,1,MPI_COMM_WORLD, &status);

    MPI_Ssend(f[3],N*N*N,MPI_DOUBLE,rank-1,1,MPI_COMM_WORLD); 
  }
  ---------------------------------
\end{verbatim}

To make a rough estimate of the total speedup let $n$ be the number of processes and $np$ be the number of cores used (assume a fixed number of cores per node). Then the speedup can be described by, to leading order,
\begin{equation}
\frac{T_{\textrm{SERIAL}}}{T_{\textrm{PARALLEL}}} = \frac{ C M N^6 T_{\textrm{FLOP}} }{4nN^3 T_{\textrm{MEM}} + C M N^6 T_{\textrm{FLOP}}/np}, 
\end{equation}
where $T_{\textrm{MEM}}$ is the average time to transfer a double precision value and $T_{\textrm{FLOP}}$ is the time for a single floating point operation. As $n$ becomes large with $N$ and $M$ fixed, more and more data transfer is required, however one would need $4n^2 p T_{MEM}/T_{FLOP} \approx C M N^3$ before the memory access time would begin to dominate the collisional computations.

We remark that this parallelization approach can work for other collisional models that take the form of a weighted convolution in Fourier space, for example collisional models with anisotropic scattering or the Landau equation \cite{GamHaaRGDAnIso, GamHaaAnIso}.


We test the this method with the sudden change in wall temperature example suggested by Aoki et al. \cite{AokiSoneNisSug91}. These authors computed this example using the BGK approximation of the collision operator and finite differences.  This was later computed in \cite{GamTha10} with a preliminary first order serial version of this conservative spectral method for the fullt nonlinear Boltzmann operator. In this problem the gas is assumed to be initially at equilibrium, and the temperature of the wall at the boundary of the domain is instantaneously changed at $t=0$. This gives rise to a discontinuous distribution function at the wall, which propagates into the domain and eventually forms a shock. This problem is difficult to compute with a Monte Carlo method due a distribution function that is far from equilibrium near the wall and the fact that the the shock develops over a long period of time.

\begin{figure}[!htbp]\label{SuddenHeat}
\includegraphics[width=0.45\linewidth]{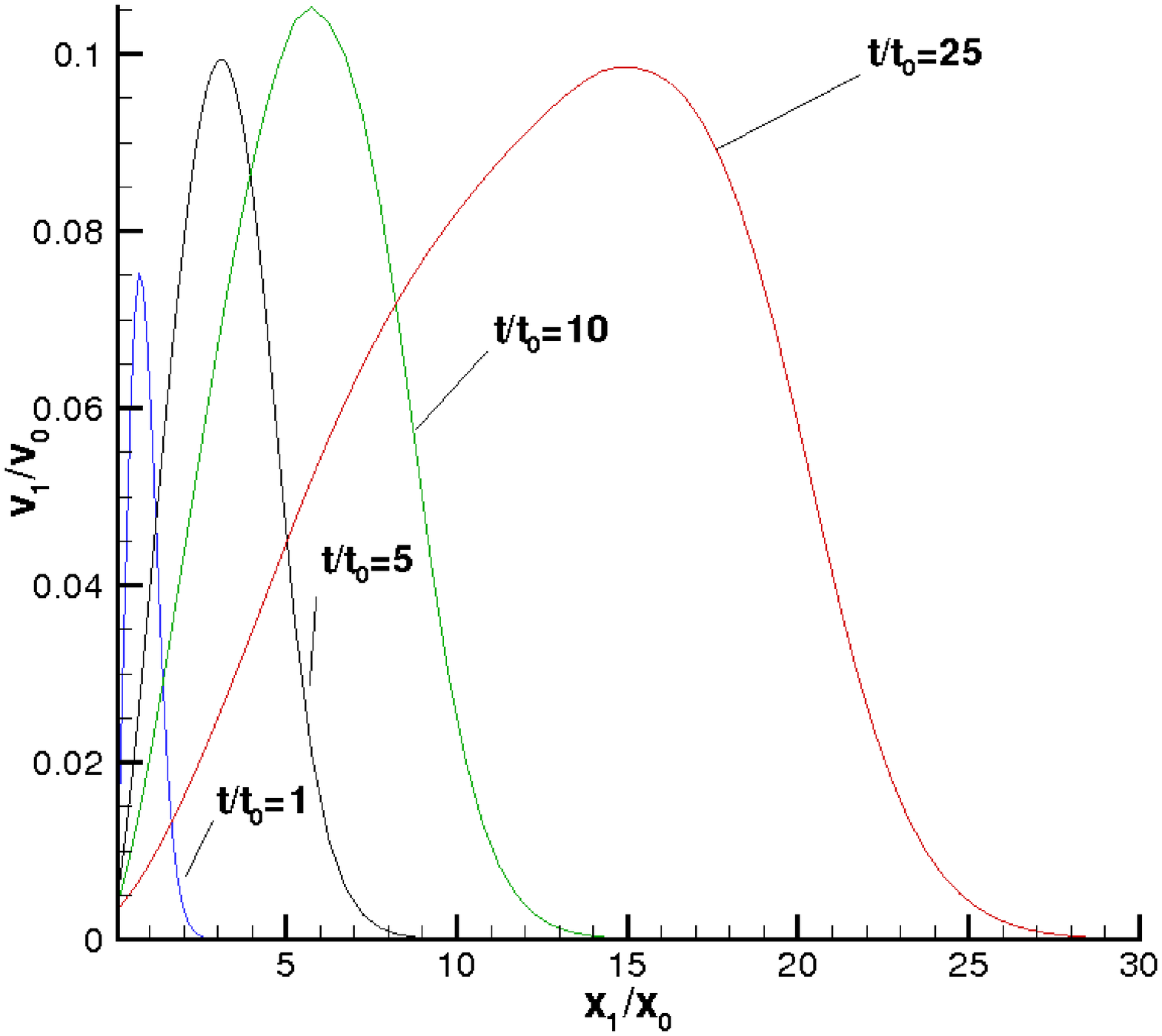} 
\includegraphics[width=0.45\linewidth]{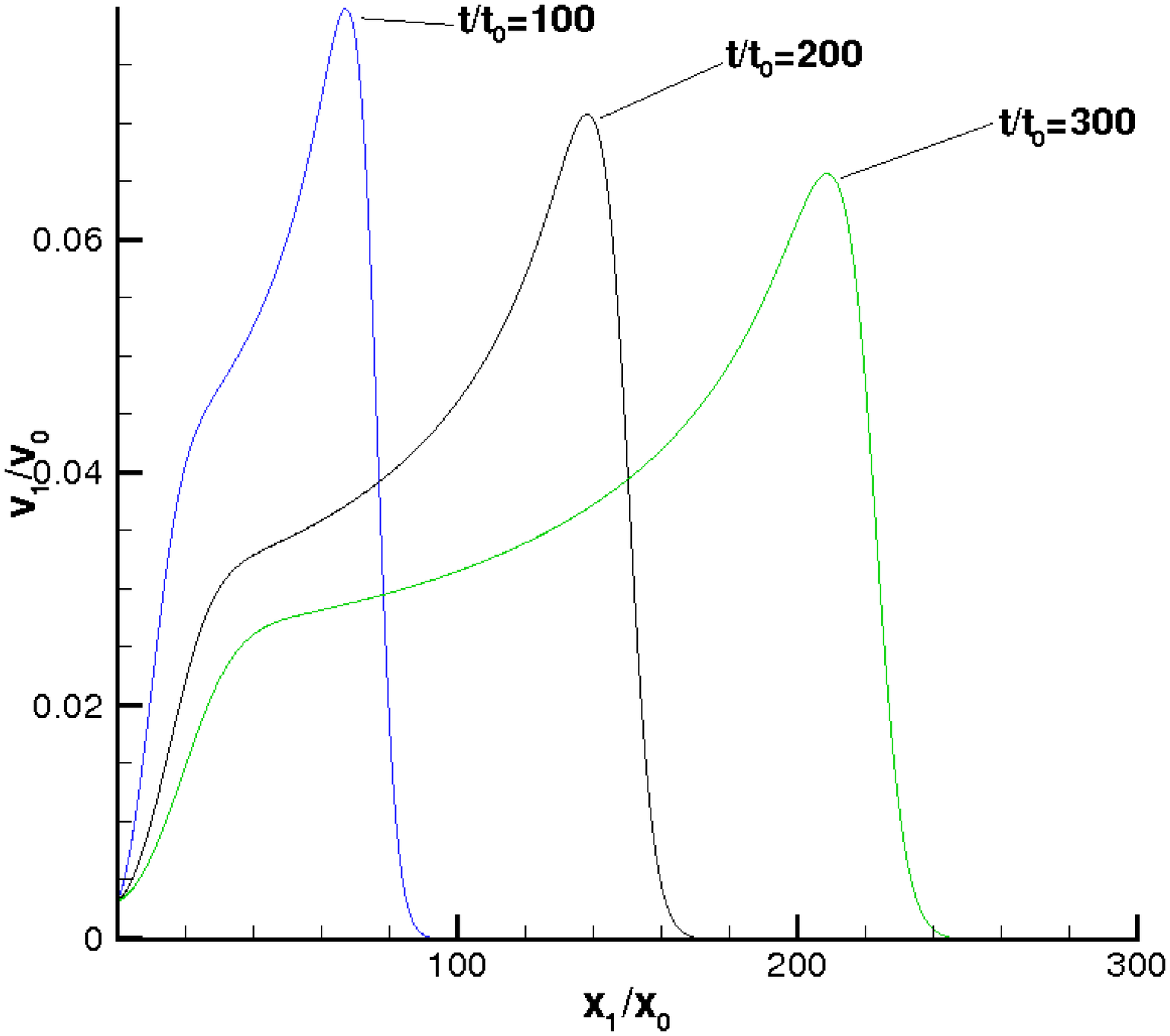} 
\caption{Formation of shock from sudden heating of wall. Evolution of bulk velocity.}
\end{figure}

In Figure \ref{SuddenHeat} we show shock formation due to a sudden change in wall temperature. Unlike the computations in \cite{AokiSoneNisSug91}, only two grid points are used per mean free path in the interior of the domain. Near the wall, the grid is refined to eight points per mean free path to better capture the finer dynamics where the distribution function is discontinuous. In both examples, we set the scaled Knudsen number to $\eps = 1$.  Note that despite the discontinuity we do not observe any Gibbs phenomena in the solution. We hypothesize that this is due to the mixing effects of the convolution weights; this will be explored in future work.

\begin{figure}[!htbp]\label{SuddenHeatTemp}
\includegraphics[width=0.45\linewidth]{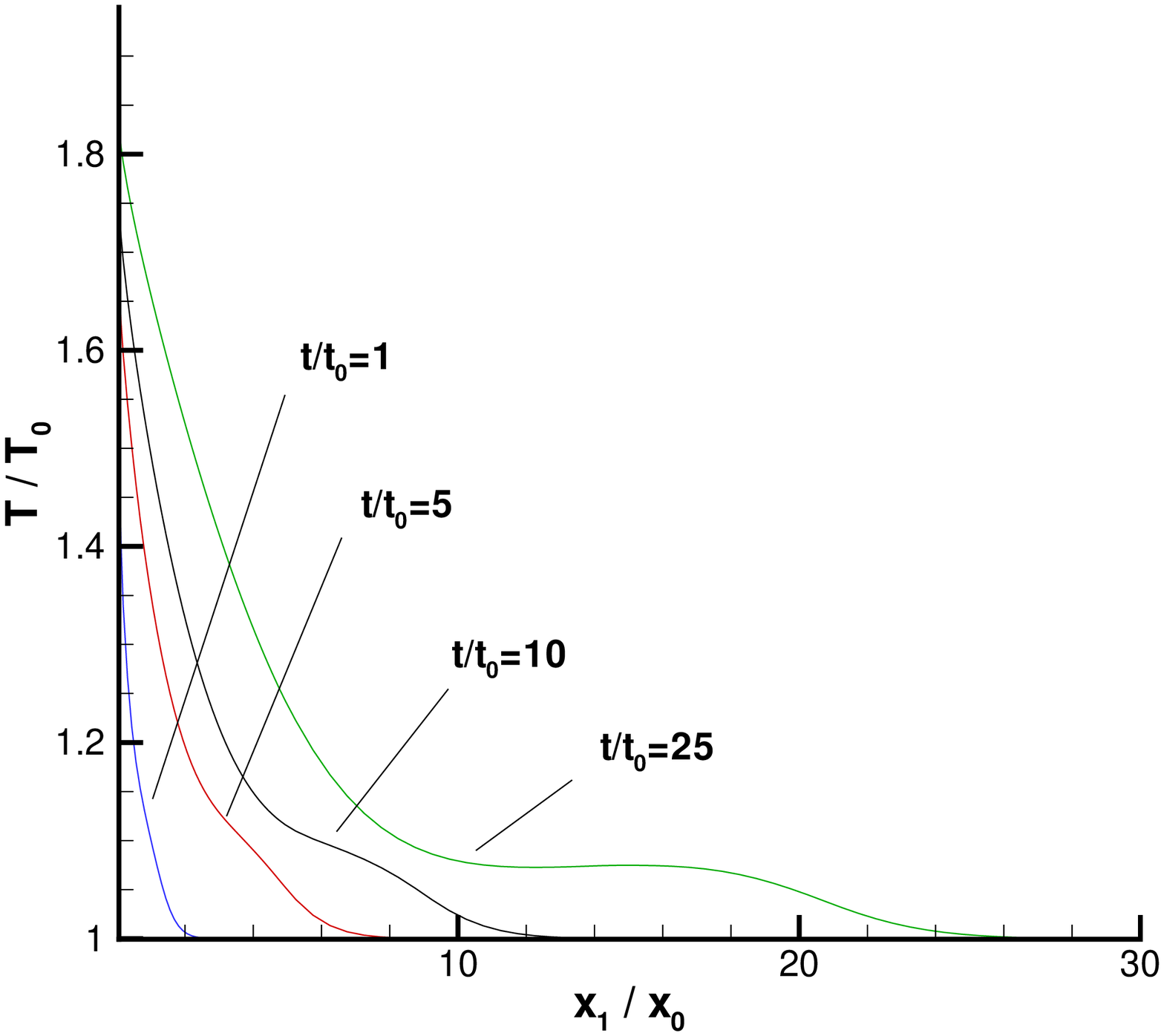} 
\includegraphics[width=0.45\linewidth]{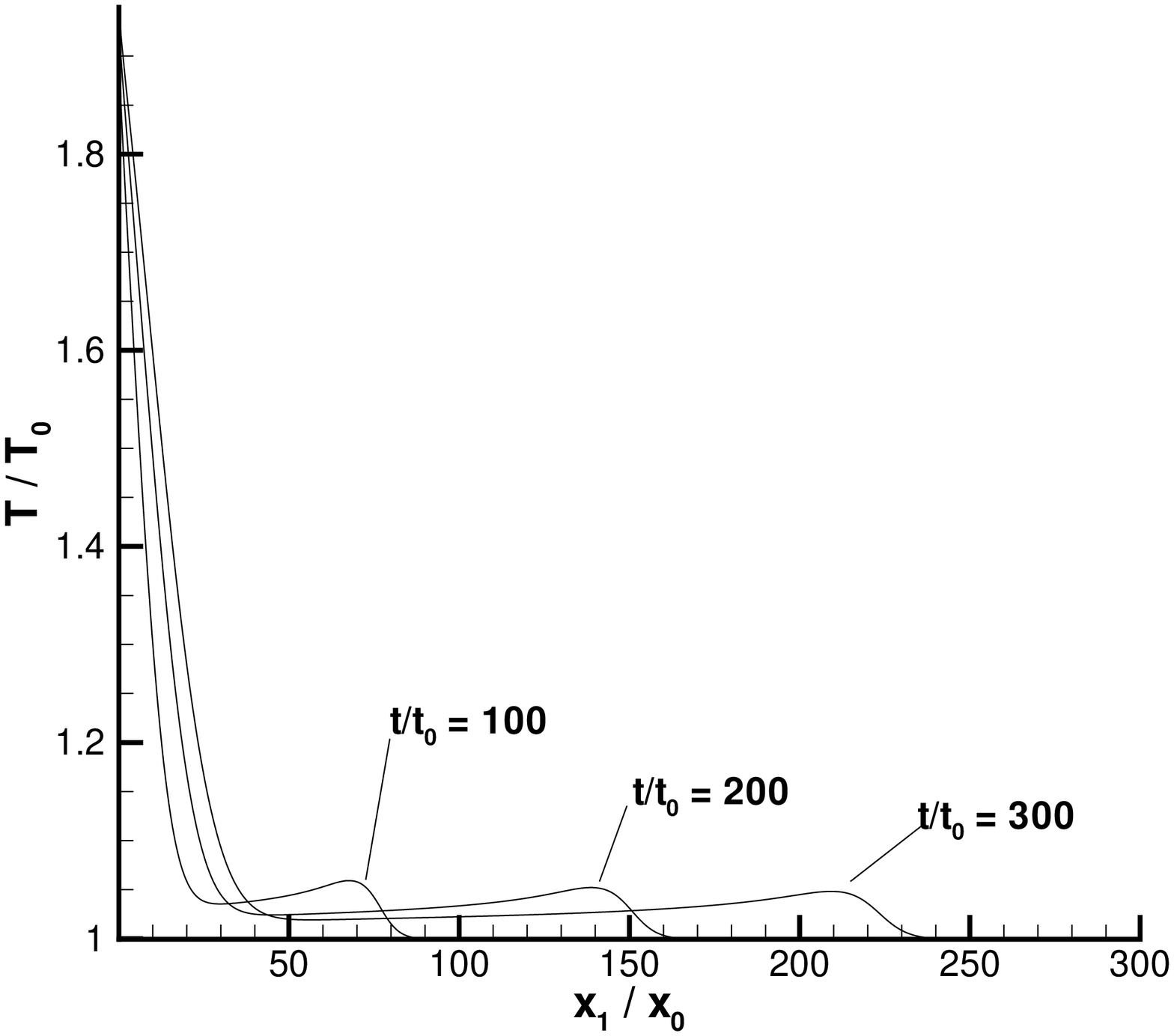} 
\caption{Formation of shock from sudden heating of wall. Evolution of kinetic temperature.}
\end{figure}

\begin{figure}[!htbp]\label{SuddenCool}
\includegraphics[width=0.45\linewidth]{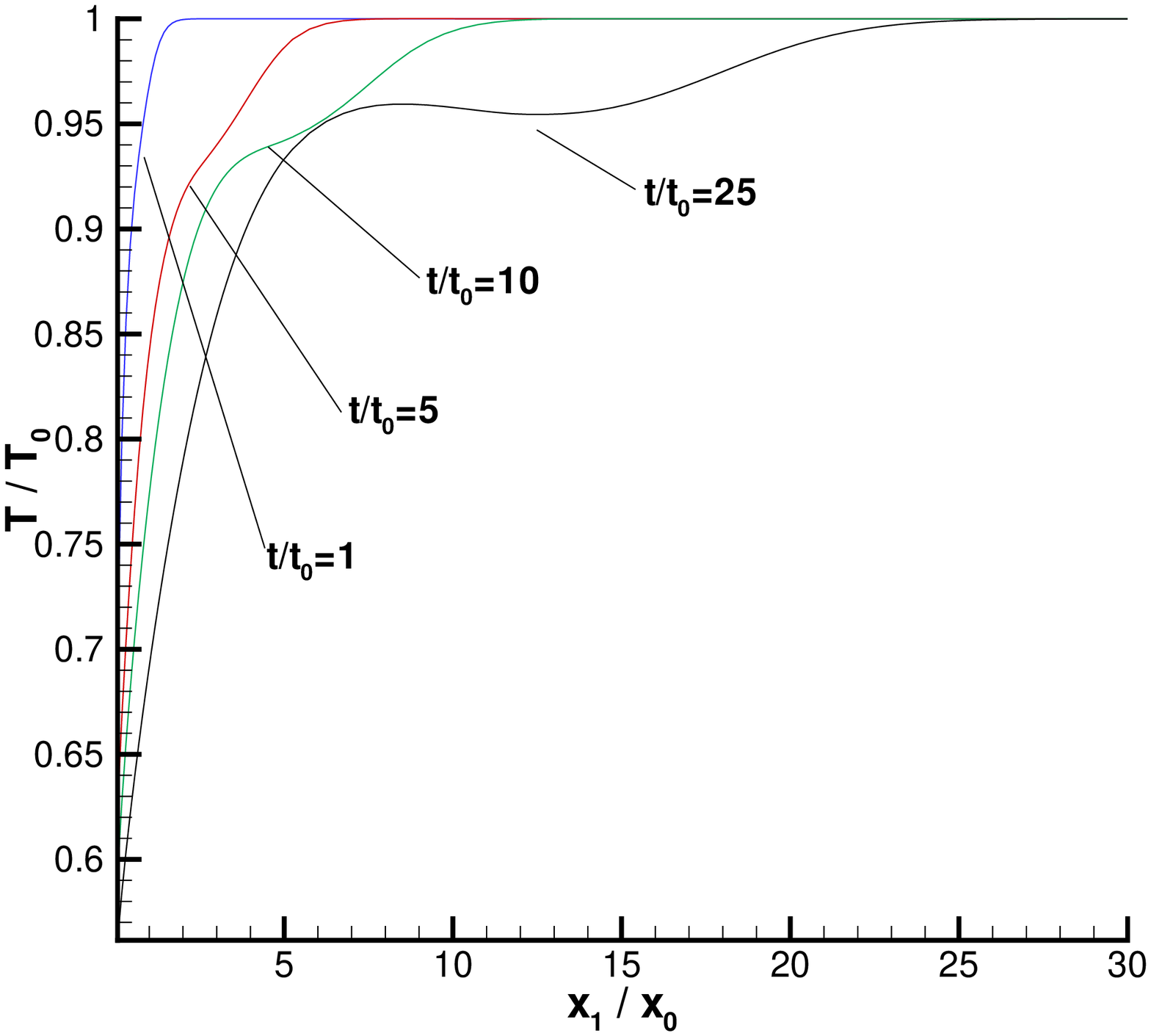} 
\includegraphics[width=0.45\linewidth]{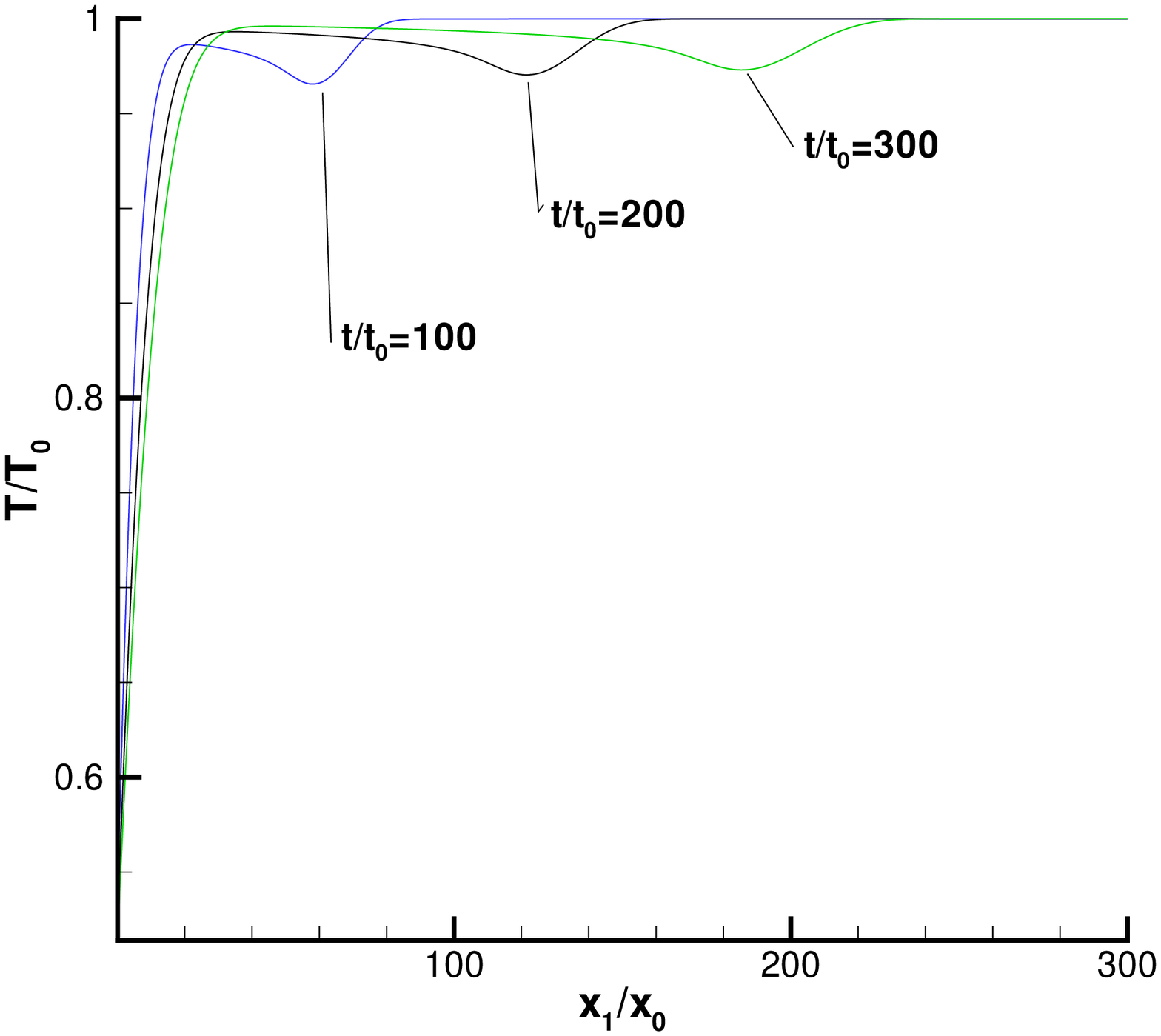} 
\caption{Results from sudden cooling of wall. Evolution of kinetic temperature.}
\end{figure}

\begin{figure}[!htbp]\label{Marginal}
\includegraphics[width=0.95\linewidth]{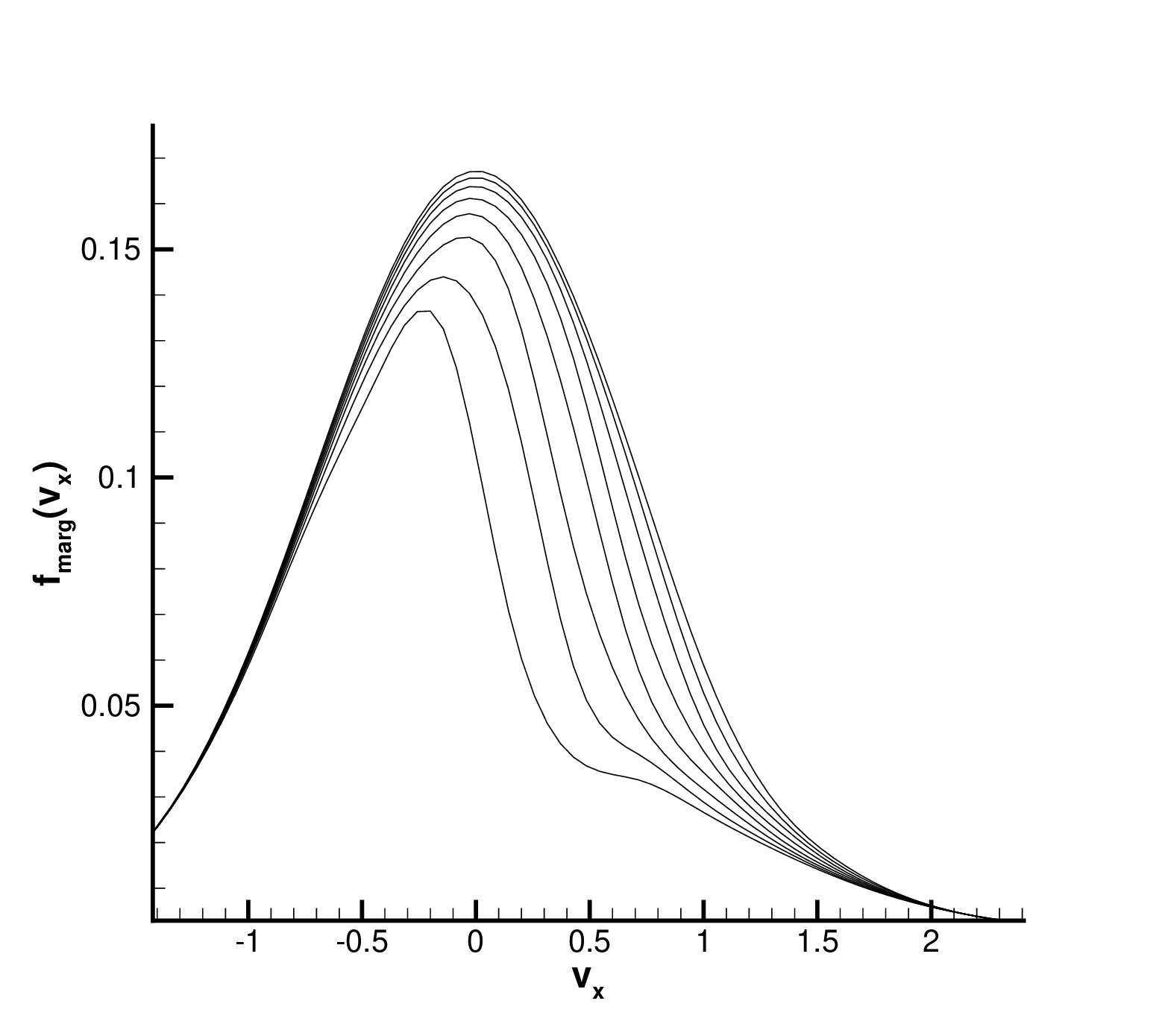} 
\caption{Sudden heating: evolution of discontinuous marginal distribution near the wall. The plots from bottom to top are the marginal distribution of $f$ in eight equispaced cells on $x=[0,1]$}
\end{figure}

In Table \ref{times} we show the wall time scaling results for the sudden heating example. We see a near perfect linear scaling as more nodes are added. Overall the openMP/MPI hybrid scheme gives a speedup of $\sim 4n$, where $n$ is the total number of nodes. Based on the openMP results above, we would expect a $8n$ increase on the Stampede processors. For the computations above, we used 32 nodes, or 512 cores, to obtain a speedup of $\sim128$ compared to the computations in \cite{GamTha10}. 

\begin{table}[!htbp]
\begin{tabularx}{\textwidth}{XXX}
\hline
nodes & cores & time (s)\\
\hline
1 & 16 & 456.313\\
2 & 32 & 235.315 \\
4 & 64 & 120.762 \\
8 & 128 & 61.345 \\
16 & 256 & 30.943 \\
32 & 512 & 15.252 \\
64 & 1024 & 7.813 \\
128 & 2048 & 4.042 \\
\hline
\end{tabularx}
\caption{Computational time for a single timestep in sudden heating example from Figure 1.}
\centering
\label{times}
\end{table}


\section{Conclusions} \label{sec:conclusion}
We have extended the spectral method of Gamba and Tharkabhushanam to a second order scheme with a nonuniform grid in physical space, and investigated its scaling to high performance computing. The method showed nearly linear speedup across nodes when applied to a large computation problem solved on a supercomputer. However, at some point memory access will still become a problem, not with transfer between nodes but with memory access on a single node. The most expensive object to store is the six dimensional weight array $\widehat{G}(\zeta,\xi)$, and if $N$ becomes too large it may not fit on a single node's memory, significantly slowing down computation. At that point, it may be faster to simply compute the weights on the fly, as flops are much cheaper than memory accesses on the large distributed systems. For the isotropic scattering cross sections considered in this paper this would require computation of at most a one dimensional integral for each velocity grid point, and in the hard sphere and Maxwell molecules we have an exact formula to find the weights. Memory management will become especially important when implementing the code on MICs, which have much more processing power but relatively smaller memory. This will be a subject of future work.

\section{Acknowledgments}
Thanks to Irene Gamba for introducing me to this method, and thanks to Andrew Christlieb for opening my eyes to parallel programming. This work has been supported by the NSF under grant number DMS-0636586. 

\bibliographystyle{siam}   
\bibliography{Boltz}
\end{document}